# CRITICAL RANDOM HYPERGRAPHS: THE EMERGENCE OF A GIANT SET OF IDENTIFIABLE VERTICES

By Christina Goldschmidt

*Université Pierre et Marie Curie (Paris 6)*

We consider a model for random hypergraphs with *identifiability*, an analogue of connectedness. This model has a phase transition in the proportion of identifiable vertices when the underlying random graph becomes critical. The phase transition takes various forms, depending on the values of the parameters controlling the different types of hyperedges. It may be continuous as in a random graph. (In fact, when there are no higher-order edges, it is exactly the emergence of the giant component.) In this case, there is a sequence of possible sizes of "components" (including but not restricted to $N^{2/3}$). Alternatively, the phase transition may be discontinuous. We are particularly interested in the nature of the discontinuous phase transition and are able to exhibit precise asymptotics. Our method extends a result of Aldous [*Ann. Probab.* **25** (1997) 812–854] on component sizes in a random graph.

**1. Poisson random hypergraphs.** The emergence of the giant component in a random graph is now a well-understood phenomenon (see [5]). The purpose of this paper is to demonstrate that an analogous, but richer, phenomenon occurs in random hypergraphs. We employ stochastic process methods of the type described in [6].

We use the framework of Poisson random hypergraphs introduced by Darling and Norris [3]. Suppose we are given a set of vertices $V$ of size $N$. Denote the power set of $V$ by $\mathcal{P}(V)$. We define a Poisson random hypergraph with parameters $(\beta_k : k \geq 2)$ by a random map $\Lambda : \mathcal{P}(V) \to \mathbb{Z}^+$ such that

$$\Lambda(A) \sim \text{Poisson}\left(N\beta_k \Big/ \binom{N}{k}\right) \qquad \text{whenever } |A| = k.$$

$\Lambda(A)$ is the number of hyperedges of size $k$ (or "$k$-edges") over the set $A$. The numbers of hyperedges over different subsets of the vertex set are independent. We allow multiple edges and the distribution of $\Lambda(A)$ depends only on $|A|$. Define a generating function, $\beta(t) = \sum_{k=2}^{\infty} \beta_k t^k$. Throughout this paper, we shall assume that $\beta'(1) < \infty$, which ensures that each vertex is contained in a finite number of edges.

We now proceed by defining identifiability. This analogue of connectedness was first introduced by Darling and Norris [3] and further studied by Darling, Levin and Norris [2]. We put 1-edges (or *patches*) on certain arbitrarily chosen vertices, which make those vertices act like roots of "components." We define the set of identifiable vertices to be the smallest set satisfying the following recursive conditions. First, vertices with patches on them are identifiable. Furthermore, a vertex without a patch is identifiable if there exists $r \geq 2$ such that it is contained a hyperedge of size $r$ and the other $r-1$ vertices are all identifiable. In a random graph (with patches), the identifiable vertices are those which are in the same component as a vertex with a patch on it. Thus, the patches "pick out" some of the components.

It is useful to have an algorithm for finding the identifiable vertices in a hypergraph. Pick a patch at random and delete it and the vertex $v$ underneath it. Collapse all of the other hyperedges over $v$ down onto their remaining vertices so that, for example, a 4-edge over $\{u, v, w, x\}$ becomes a 3-edge over $\{u, w, x\}$. In particular, any 2-edges including $v$ become new patches. Repeat until there are no patches left. The set of identifiable vertices consists of those vertices deleted. It turns out that the order of deleting patches in this collapse procedure does not affect the ultimate set of identifiable vertices (see [3]).

We will also consider the vertices identifiable from a particular vertex. Define the *domain* of $v$ to be the set of vertices identifiable when a single patch is put on the hypergraph, at $v$. An equivalent definition is as follows: $w$ is in the domain of $v$ if and only if either (a) $w = v$ or (b) there exists a sequence $v_0 = v, v_1, v_2, \ldots, v_r = w$ for some $r \geq 1$ such that for each $1 \leq i \leq r$ there exists a hyperedge consisting of $v_i$ and some nonempty subset of $\{v_0, v_1, \ldots, v_{i-1}\}$. In a graph, the domain of $v$ is the same as its component and, indeed, domains will play the role of components in what follows. Note, however, that in a general hypergraph identifiability from a vertex is not a symmetric property. It is perfectly possible for $w$ to be in the domain of $v$ without $v$ being in the domain of $w$. Because of this lack of symmetry, the analogy with the concept of a component is incomplete. We observe, nonetheless, that in a Poisson random hypergraph with $\beta_2 > 0$ there is an underlying random graph and the domains of any two vertices in the same underlying 2-edge component are the same.

THEOREM 1.1 ([2]). *Let $D_N$ be the size of the domain of an arbitrarily chosen vertex in a Poisson random hypergraph on $N$ vertices with $\beta_2 > 0$.*

*Define* $t^* = \inf\{t \geq 0 : \beta'(t) + \log(1-t) < 0\}$ *and suppose that there are no zeros of* $\beta'(t) + \log(1-t)$ *in* $(0, t^*)$. *Then*

$$\frac{1}{N} D_N \xrightarrow{d} t^* \mathbb{1}_{\{M=\infty\}}$$

*as* $N \to \infty$, *where* $M$ *has the* Borel$(2\beta_2)$ *distribution* [*i.e., the distribution of the total population of a branching process with* Poisson$(2\beta_2)$ *offspring*].

Let us first discuss the meaning of this result for a random graph. The critical value for the emergence of a giant component is $\beta_2 = 1/2$. The significance of $M$ is that the component of the specific vertex on which we put our patch has a size which converges in distribution to $M$ as $N \to \infty$ (this reflects the fact that small components look very much like branching process family trees). If $\beta_2 \leq 1/2$, the largest component is of size $o(N)$. Thus, putting one patch onto the graph identifies at most $o(N)$ vertices and so $D_N/N \xrightarrow{p} 0$. This corresponds to the fact that for $\beta_2 \leq 1/2$, $M < \infty$ almost surely. Consider now the case $\beta_2 > 1/2$. Either the patch falls on the giant component [with probability $\mathbb{P}(M = \infty)$] and identifies a positive proportion, $t^*$, of the vertices, or it falls on a small component [with probability $\mathbb{P}(M < \infty)$] and identifies only $o(N)$. The theorem tells us that this limiting justification works even in the presence of higher-order edges but that the precise proportion identified depends on the parameters of those higher-order edges. Thus, Theorem 1.1 characterizes a phase transition in the proportion of identifiable vertices for a random hypergraph with a single patch.

For a random graph, $t^*$ represents the proportion of vertices in the giant component (and note that $t^* = 0$ for $\beta_2 \leq 1/2$). In a random hypergraph, $t^*$ represents the proportion of vertices in a "giant domain." Note that it is not clear that there is a unique such domain (although it is clear from Theorem 1.1 that any such domain must contain the vertices of any giant 2-edge component). However, it seems that a giant domain is close to being unique in that all giant domains contain an asymptotic proportion $t^*$ of the vertices. In a random hypergraph, we have $t^* = 0$ for $\beta_2 < 1/2$ but we may have $t^* > 0$ for $\beta_2 = 1/2$. To be precise, if $\beta_2 = 1/2$ and $\beta_3 > 1/6$, then $t^* > 0$, whereas if $\beta_2 = 1/2$ and $\beta_3 < 1/6$, then $t^* = 0$. If $\beta_3 = 1/6$, we must look at whether $\beta_4$ is less than or greater than $1/12$. In general, there exists a sequence of "critical" values for the $\beta_j$'s such that if there exists $k$ such that for $2 \leq j \leq k-1$, $\beta_j = 1/j(j-1)$ and $\beta_k > 1/k(k-1)$, then $t^* > 0$, whereas if there exists $k$ such that for $2 \leq j \leq k-1$, $\beta_j = 1/j(j-1)$ and $\beta_k < 1/k(k-1)$, then $t^* = 0$. [Note that the case $\beta_j = 1/j(j-1)$ for all $j \geq 2$ is explicitly excluded by the assumption $\beta'(1) < \infty$.] So it appears that, in some sense, a giant domain may already be present at the critical point in a random hypergraph (although we have probability 0 of hitting it with our single patch). Thus, the random hypergraph phase transition can be discontinuous, in that $t^*$ may not be a continuous function of $\beta_2$ at $\beta_2 = 1/2$, whereas the random graph phase transition is always continuous.

In order to investigate the random hypergraph phase transition further, we will consider what happens when, instead of a single patch, we put $\omega(N)$ patches on the critical hypergraph, where $\omega(N)/N \to 0$ as $N \to \infty$. We will add a patch to the hypergraph, collapse as far as possible and then add another patch on a vertex chosen uniformly at random from those remaining whenever needed to keep going. Is there ever a function $\omega(N)$ such that we identify $\Theta(N)$ vertices (i.e., a giant set of vertices)?

## 2. Results.

Let $\alpha(k) = (2k-4)/(2k-3)$ and

$$W^k(t) = B(t) + \frac{1}{k-1}(k(k-1)\beta_k - 1)t^{k-1},$$

where $(B(t))_{t\geq 0}$ is a standard Brownian motion.

THEOREM 2.1. *Consider a Poisson random hypergraph on $N$ vertices.*

(i) *Suppose that there exists $k \geq 3$ such that for $2 \leq j \leq k-1$, $\beta_j = 1/j(j-1)$ and $\beta_k < 1/k(k-1)$. Let $X_N$ be the number of vertices identified when $\omega(N)$ patches are added to the hypergraph one by one, as necessary, where $\omega(N)/N \to 0$ as $N \to \infty$. Then we have*

$$\frac{1}{N} X_N \xrightarrow{p} 0 \tag{2.1}$$

*as $N \to \infty$. Recall that $D_N$ is the size of the domain of a randomly chosen vertex [so that $D_N$ is the same as $X_N$ when $\omega(N) = 1$]. Then for any $\varepsilon > 0$, there exists $C$ such that for all sufficiently large $N$,*

$$\mathbb{P}(N^{-\alpha(k)} D_N < C) \geq 1 - \varepsilon. \tag{2.2}$$

(ii) *Suppose now that there exists $k \geq 3$ such that for $2 \leq j \leq k-1$, $\beta_j = 1/j(j-1)$ and $\beta_k > 1/k(k-1)$. Let $A_N^\delta$ be the number of patches we need to add one by one, as necessary, until we have identified more than $N\delta$ vertices, for $\delta > 0$. Then for all $\delta > 0$ sufficiently small,*

$$N^{-\alpha(k)/2} A_N^\delta \xrightarrow{d} -\inf_{t \geq 0} W^k(t). \tag{2.3}$$

*Let $X_N$ be the total number of vertices identified when we add patches one by one as before until at least $N\delta$ vertices have been identified (i.e., the number of vertices identified when $A_N^\delta$ patches are added), for $\delta$ sufficiently small that (2.3) is satisfied. Then*

$$\frac{1}{N} X_N \xrightarrow{p} t^* \tag{2.4}$$

*as $N \to \infty$, where $t^* = \inf\{t \geq 0 : \beta'(t) + \log(1-t) < 0\}$.*

**3. Breadth-first walk.** In order to track the process of collapse (adding patches whenever they are needed to keep going), we construct an extension to hypergraphs of the usual breadth-first walk on graphs. (Note that this is a different extension from that used in [7].)

Consider any hypergraph on $N$ vertices, with no patches and an arbitrary numbering of the vertices. Then we may define the *breadth-first ordering* as follows:

1. Take the lowest-numbered vertex, call it $v(1)$ and put a patch on it. Define the *children* of vertex $v(1)$ to be those vertices connected to it by a 2-edge. Suppose that $v(1)$ has $c(1)$ children. Number them $v(2), \dots, v(c(1)+1)$, retaining the ordering of the original labels. Now collapse the patch on $v(1)$, leaving patches on all of its children and any higher-order edges collapsed onto their remaining vertices.
2. Now look at $v(2)$. Label its children as $v(c(1)+2), \dots, v(c(1)+c(2)+1)$, where, in general, we define the children of a vertex to be those vertices connected to it by a 2-edge which have not yet been renumbered. Note that some of these children may only just have appeared as a result of the collapse of vertex $v(1)$. For example, in Figure 1, $v(3)$ is the child of $v(2)$ but only becomes visible as such after the deletion of $v(1)$.
3. Continue in this way, collapsing the vertices in numerical order [so the next one to consider is $v(3)$]. When we run out of patches, pick the next lowest-numbered vertex in the old ordering, put a patch on it and proceed as before. The process terminates when there are no more vertices to consider.

So, loosely speaking, we number within levels of an underlying tree before moving further from the "root," the only complication being that the children of $v(i)$ may only all be visible after the deletion of vertex $v(i-1)$.

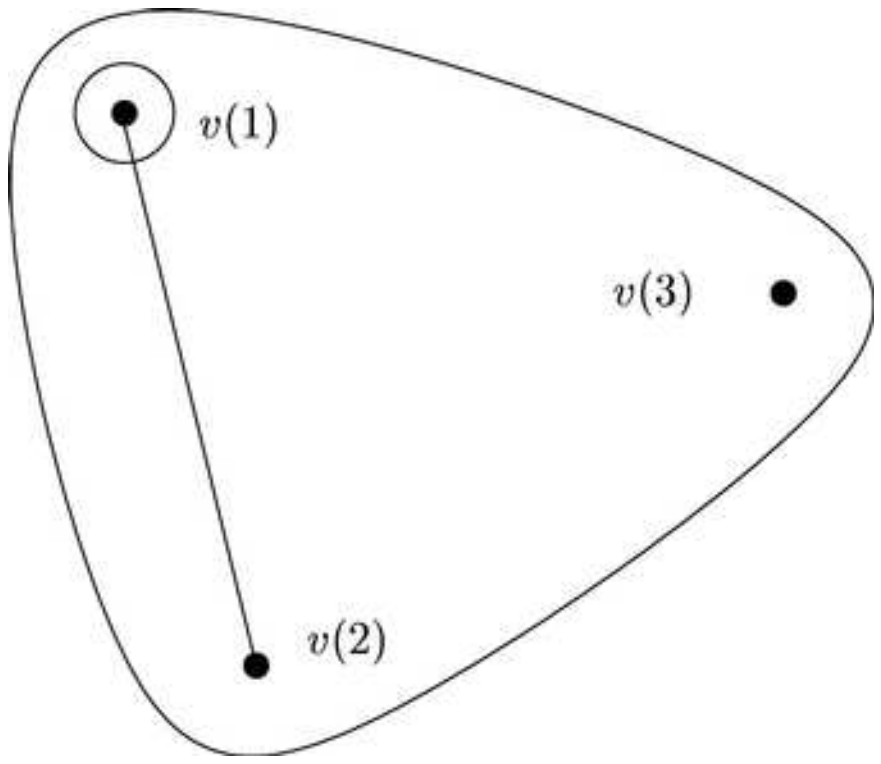

FIG. 1. *Children of vertices can appear during the process of collapse.*

Now we can define a walk $(z(i))_{0\leq i\leq N}$ on the integers associated with this hypergraph by

$$z(0)=0,$$
$$z(i)=z(i-1)+c(i)-1, \qquad i\geq 1,$$

where, as before, $c(i)$ is the number of children of vertex $v(i)$ in the breadth-first ordering. Then $i$ is the number of vertex deletions (we will also refer to this as a *time*) and $z(i)$ is the number of patched vertices on the hypergraph after the $i$th vertex deletion, minus the number of patches added before the $(i+1)$st deletion to keep the process going. The process $(z(i))_{0\leq i\leq N}$ is called the breadth-first walk.

Thus, for a *random* hypergraph on $N$ vertices, we obtain a random walk, $(Z_N(i))_{i\geq 0}$, on the integers which summarizes information about the hypergraph. [In the sequel, we will refer to the random number of children of vertex $v(i)$ as $C_N(i)$.] Most importantly, the number of vertices which are identifiable from the patches we add are coded as excursions above past minima in the breadth-first walk. This is because the breadth-first walk picks out an underlying forest structure, with each tree naturally having one more vertex than it has edges. For more details, see [1].

It will be useful later to have some notation for the number of patches added to the hypergraph so far. Let $P_N(0)=1$ and, for $i\geq 1$,

$$P_N(i)=1-\min_{j\leq i} Z_N(j).$$

Then $P_N(i)$ is the number of patches added before the $(i+1)$st vertex deletion. $Z_N(i)+P_N(i)$ is the actual number of patches on the hypergraph just after the deletion of $v(i)$ and is always greater than or equal to 1. Thus, we have the alternative representation

$$P_N(i)=1+\sum_{i=1}^{N-1}\mathbb{1}_{\{Z_N(i-1)+P_N(i-1)=1,\ C_N(i)=0\}}. \tag{3.1}$$

Recall that $\alpha(k)=(2k-4)/(2k-3)$ and that

$$W^k(t)=B(t)+\frac{1}{k-1}(k(k-1)\beta_k-1)t^{k-1},$$

where $(B(t))_{t\geq 0}$ is a standard Brownian motion. Then our key result is the following:

THEOREM 3.1. *Suppose that $(Z_N(i))_{0\leq i\leq N}$ is the breadth-first walk on the Poisson random hypergraph on $N$ vertices and that there exists a $k\geq 3$ such that $\beta_j=1/j(j-1)$ for all $2\leq j\leq k-1$. Rescale by defining*

$$\bar{Z}_N(t)=N^{-\alpha(k)/2}Z_N(\lfloor N^{\alpha(k)}t\rfloor).$$

*Then*

$$\bar{Z}_N \xrightarrow{d} W^k$$

*as* $N \to \infty$ *in* $D[0, \infty)$.

Note that here the convergence is uniform on compact time-intervals. The proof of this result is deferred to Section 6 to enable us to first interpret its implications for the hypergraph.

**4. Consequences.** If there exists $k$ such that $\beta_j = 1/j(j-1)$ for all $2 \le j \le k-1$ and $\beta_k < 1/k(k-1)$, then the limit process has a negative drift which increases in magnitude with time. Thus, the process keeps hitting its previous minima, on average resulting in smaller and smaller numbers of identifiable vertices per patch added. This is very like what we see in Theorem 3 of [1] where the components of a random graph appear in size-biased order. In the critical random graph case, $\beta_3 = 0 < 1/6$ and the components are of size $\mathcal{O}(N^{2/3})$, as is well known. However, in the random hypergraph, there is a whole series of critical scalings ($N^{2/3}$, $N^{4/5}$, $N^{6/7}$, ...) which can be attained by suitable adjustments of the parameters $\beta_3, \beta_4, \ldots$. Thus, the random hypergraph demonstrates much richer behavior than the random graph.

If there exists $k$ such that $\beta_j = 1/j(j-1)$ for all $2 \le j \le k-1$ and $\beta_k > 1/k(k-1)$, then the process $W^k$ has positive drift and so there is a (random) last time that it hits its own minimum. This signals the start of a giant excursion which is too big to be measured on the scale of $N^{\alpha(k)}$. We wish to prove that the domain which this excursion represents is, in fact, of size comparable to $N$. In order to do this, we will show that the giant excursion has length at least $N\delta$ for all sufficiently small $\delta > 0$. This will then allow us to prove a fluid limit theorem for the breadth-first walk; the length of the excursion of the fluid limit above 0 gives us the asymptotic size of the giant set of identifiable vertices. We will also discuss the fluctuations of the breadth-first walk around this fluid limit.

**5. The giant set of identifiable vertices.** For ease of notation, define $\mu_k = k(k-1)\beta_k - 1$. We will now fix $k \ge 3$ and look at the case $\mu_k > 0$, $\mu_j = 0$ for $2 \le j \le k-1$ in more detail. First, we state a proposition which will be of use to us later:

PROPOSITION 5.1. *For* $W^k$ *defined as in Section* 2,

$$\mathbb{P}(W^k(R^2) > R) \to 1$$

*as* $R \to \infty$.

PROOF. We have

$$\begin{aligned}\mathbb{P}(W^k(R^2) > R) &= \mathbb{P}\left(B(R^2) + \frac{\mu_k}{k-1} R^{2(k-1)} > R\right) \\ &= 1 - \Phi\left(1 - \frac{\mu_k}{k-1} R^{2k-3}\right) \\ &\to 1\end{aligned}$$

as $R \to \infty$, where $\Phi$ is the standard Normal distribution function. □

Thus, it is a corollary of Theorem 3.1 that the event $\{Z_N(\lceil N^{\alpha(k)} R^2 \rceil) > RN^{\alpha(k)/2}\}$ has asymptotically high probability in the sense that

$$\lim_{R\to\infty} \lim_{N\to\infty} \mathbb{P}(Z_N(\lceil N^{\alpha(k)} R^2 \rceil) > RN^{\alpha(k)/2}) = 1. \tag{5.1}$$

Recall that $P_N(i) = 1 - \min_{j \le i} Z_N(i)$ is the number of patches added before the $(i+1)$st deletion to keep the process going. Then, by Theorem 3.1 and the continuous mapping theorem (Corollary 3.1.9 of [4]), $N^{-\alpha(k)/2} P_N(\lfloor N^{\alpha(k)} t \rfloor) \xrightarrow{d} -\inf_{0 \le s \le t} W^k(s)$. Because of the positive drift of $W^k$, we have

$$\mathbb{P}\left(\inf_{s \ge 0} W^k(s) < -R\right) \to 0$$

as $R \to \infty$ and so it is clear that

$$\lim_{R\to\infty} \lim_{N\to\infty} \mathbb{P}(P_N(\lceil N^{\alpha(k)} t \rceil) < RN^{\alpha(k)/2}) = 1 \tag{5.2}$$

for any value of $t$.

Define

$$S_N^R = \inf\{i \ge N^{\alpha(k)} R^2 : Z_N(i) \le 0\}.$$

THEOREM 5.2. *There exists $\delta > 0$ such that*

$$\lim_{R\to\infty} \lim_{N\to\infty} \mathbb{P}(S_N^R \le N\delta) = 0.$$

Essentially, by time $\lceil N^{\alpha(k)} R^2 \rceil$, the process $Z_N$ is, with high probability, in an excursion above its last minimum of length $\Theta(N)$.

Define

$$z(t) = 1 - t - \exp(-\beta'(t)) \tag{5.3}$$

and recall that $t^* = \inf\{t \ge 0 : z(t) < 0\}$. Assume that there are no zeros of the function $z(t)$ in $(0, t^*)$ and note that $t^* < 1$. Let $\tilde{Z}_N$ be a modified version of $Z_N$ such that no more patches are added after time $\lfloor N\delta \rfloor$, for $\delta$ as in Theorem 5.2. Thus, $\tilde{P}_N(\lfloor Nt \rfloor) = P_N(\lfloor N\delta \rfloor)$ for all $t \ge \delta$ and the first time that $\tilde{Z}_N$ goes below its past-minimum after time $\lfloor N\delta \rfloor$, it stops evolving and stays constant. Let $\tilde{z}(t) = z(t \wedge t^*)$. Theorem 5.2 allows us to prove a fluid limit for this modified version of $Z_N$.

Theorem 5.3. *For all* $\varepsilon > 0$,

$$\lim_{N\to\infty} \mathbb{P}\left( \sup_{0\le t\le 1} \left| \frac{1}{N} \tilde{Z}_N(\lfloor Nt \rfloor) - \tilde{z}(t) \right| > \varepsilon \right) = 0.$$

This implies that for any $0 < \sigma < t^*$,

$$\sup_{0\le t\le \sigma} \left| \frac{1}{N} Z_N(\lfloor Nt \rfloor) - z(t) \right| \xrightarrow{p} 0.$$

In addition to Theorem 5.3, we have a functional central limit theorem, which describes the fluctuations of the breadth-first walk around its fluid limit.

Theorem 5.4. *For any* $0 < \sigma < t^*$,

$$\frac{1}{\sqrt{N}} (Z_N(Nt) - Nz(t))_{0\le t\le \sigma} \xrightarrow{d} (X_t)_{0\le t\le\sigma},$$

*where*

$$X_t = \exp(-\beta'(t)) \int_0^t \exp(\beta'(s))\, dG_s$$

*and* $(G_t)_{t\ge 0}$ *is a Gaussian process such that if* $(B_t)_{t\ge 0}$ *is a standard Brownian motion, then* $(G_t)_{t\ge 0} \sim (B_{z(t)+t})_{t\ge 0}$.

Remark. This result is consistent with Theorem 3.1; the scaling there gives a zoomed-in version.

Assuming Theorems 3.1, 5.2 and 5.3 (which are proved in Section 6), we may now prove Theorem 2.1.

Proof of Theorem 2.1. (i) Darling and Norris [3] studied the limiting proportion of identifiable vertices in a hypergraph with a Poisson($N\beta_1$) number of patches, where $\beta_1 > 0$ and the patches are placed on vertices chosen uniformly at random (with replacement) right at the start, before any collapse has occurred. In their Theorem 2.1, they show that this limiting proportion is

$$t^*_{\beta_1} = \inf\{t \ge 0 : \beta_1 + \beta'(t) + \log(1-t) < 0\}.$$

Note that if there exists a $k \ge 3$ such that $\beta_j = 1/j(j-1)$ for $2 \le j \le k-1$ and $\beta_k < 1/k(k-1)$, then $t^*_{\beta_1} \to 0$ as $\beta_1 \to 0$. We will exploit this result to show (2.1).

The first thing we need to do is to find a way to compare the situation where we put all our patches on the hypergraph right at the start with the situation where we put them on one by one when needed to continue the process of collapse. We can couple the two situations as follows. Fix a particular realization of a Poisson random hypergraph on $N$ vertices with parameters all critical up to a certain point and then one subcritical. Take the breadth-first construction, as outlined in Section 3. Instead of putting the patches on root vertices which are always the next-lowest numbered when we come to the end of a domain, we try to put our next patch on a vertex chosen uniformly at random from $\{1,\dots,N\}$. Of course, there is some probability that the vertex we choose has already been identified, in which case we keep trying in the same manner until we find a vertex which has not yet been identified. Clearly, this takes a geometric number of trials. Then, within a domain, we continue as before in the breadth-first order. Fix $\delta > 0$ and stop the process of identification if ever we have reached the end of a domain and more than $N\delta$ vertices have been identified.

Suppose we identify the domains of up to $\omega(N)$ root vertices (before having identified $N\delta$ vertices). Let $\pi_N$ be the number of vertices (possibly counting some more than once) on which we try to put patches, including the up to $\omega(N)$ successful placings. Each of these $\pi_N$ vertices is drawn uniformly at random from $\{1,\dots,N\}$ and putting $\pi_N$ patches down on them right at the start would have identified the same vertices as putting the up to $\omega(N)$ patches on one by one when needed. Then taking $G_1,\dots,G_{\omega(N)}$ to be independent and identically distributed Geometric$(1-\delta)$ random variables, we have

$$\pi_N \le_{\mathrm{st}} \sum_{i=1}^{\omega(N)} G_i,$$

because the proportion of vertices already identified each time we try to find a root vertex is always less than $\delta$.

Let $\pi_N^\varepsilon$ be an independent Poisson$(N\varepsilon)$ random variable and let $X_N^\varepsilon$ be the number of vertices identified when $\pi_N^\varepsilon$ patches are placed on the hypergraph right at the start. Then, for any $\delta > 0$,

$$\begin{aligned}\mathbb{P}(X_N > N\delta) &\le \mathbb{P}(X_N > N\delta | \pi_N < \pi_N^\varepsilon) + \mathbb{P}(\pi_N \ge \pi_N^\varepsilon)\\ &\le \mathbb{P}(X_N^\varepsilon > N\delta) + \mathbb{P}(\pi_N \ge \pi_N^\varepsilon)\\ &\le \mathbb{P}(X_N^\varepsilon > N\delta) + \mathbb{P}\left(\sum_{i=1}^{\omega(N)} G_i \ge \pi_N^\varepsilon\right).\end{aligned}$$

The second line follows from the obvious monotonicity property that adding more patches identifies a stochastically larger number of vertices. By Theorem 2.1 of [3], we have $X_N^\varepsilon/N \xrightarrow{p} t_\varepsilon^*$ and $t_\varepsilon^* \to 0$ as $\varepsilon \to 0$. Thus, if we take $\varepsilon$ small enough that $t_\varepsilon^* < \delta$, we have that $\mathbb{P}(X_N^\varepsilon > N\delta) \to 0$ as $N \to \infty$. Moreover, as $\omega(N)/N \to 0$ as $N \to \infty$, we have that $\mathbb{P}(\sum_{i=1}^{\omega(N)} G_i \ge \pi_N^\varepsilon) \to 0$ as

$N \to \infty$, for any $\varepsilon > 0$. Thus, for any $\delta > 0$,

$$\mathbb{P}(X_N > N\delta) \to 0$$

as $N \to \infty$. This gives (2.1). Equation (2.2) follows immediately from Theorem 3.1.

(ii) $A_N^\delta = P_N(\lfloor N\delta \rfloor)$ and so Theorems 3.1 and 5.2 give (2.3). Equation (2.4) then follows from Theorem 5.3. □

## 6. Proofs.

Proof of Theorem 3.1. Our method of proof follows that of Theorem 3 of [1].

Lemma 6.1. *Let*

$$\lambda_2(N,i) = N \sum_{j=0}^{i} \beta_{j+2} \binom{i}{j} \Big/ \binom{N}{j+2}$$

*and* $\rho(N,i) = 1 - \exp(-\lambda_2(N,i))$. *Then, given* $Z_N(i-1)$ *and* $P_N(i-1)$, *the distribution of* $C_N(i)$ *is*

$$\mathrm{Bin}(N - (i-1) - Z_N(i-1) - P_N(i-1), \rho(N, i-1)).$$

Proof. Consider a Poisson random hypergraph, $\Lambda$, and delete a deterministic set $S$ of the vertices, collapsing any hyperedges with vertices in $S$ down onto their remaining vertices. Suppose $|S| = i$. Call the hypergraph on the remaining vertices $\Lambda'$. For any $A \subseteq V \setminus S$,

$$\Lambda'(A) = \sum_{B \subseteq V : A \subseteq B} \Lambda(B).$$

If $|A| = k$, this has a Poisson distribution with parameter

$$\lambda_k(N,i) = N \sum_{j=0}^{i} \beta_{k+j} \binom{i}{j} \Big/ \binom{N}{j+k}$$

and the random variables $(\Lambda'(A) : A \subseteq V \setminus S)$ inherit independence from the original hypergraph. Thus, $\Lambda'$ is another Poisson random hypergraph.

When we perform the breadth-first walk on the hypergraph we do not delete a deterministic set of vertices. We aim to show, nonetheless, that when we have deleted $i$ vertices, the probability that the number of 2-edges over the pair $\{v(i+1), w\}$, for any vertex $w$ we have not yet looked at, is $\mathrm{Poisson}(\lambda_2(N,i))$, independently for all such $w$. Start with a Poisson random hypergraph and let $\mathcal{L}$ be the set consisting of a list of the vertices. Let $\mathcal{R}$ be the set of vertices we have reached (empty to start with). Perform the following version of the breadth-first numbering algorithm:

1. Remove a vertex at random from $\mathcal{L}$, call it $v(1)$ and add it to $\mathcal{R}$.
2. Examine the number of hyperedges over the sets $\{v(1), w\}$, $w \in \mathcal{L}$. The children of $v(1)$ are those $w$ such that this number of hyperedges is 1 or greater.
3. Retaining the original ordering of the vertex-labels, call the children $v(2), \dots,$
$v(C_N(1)+1)$, where $C_N(1)$ is the number of children of $v(1)$. Remove $v(2), \dots, v(C_N(1)+1)$ from the list $\mathcal{L}$ and add them to $\mathcal{R}$.
4. Suppose we have already found the children of vertices $v(1), \dots, v(i-1)$. If there is a vertex $v(i)$ in the set $\mathcal{R}$, go to step 5. Otherwise, take a randomly chosen vertex from $\mathcal{L}$, call that $v(i)$, add it to $\mathcal{R}$ and go to step 5.
5. The children of $v(i)$ are those $w \in \mathcal{L}$ such that the number of original hyperedges over the set $\{v(i), w\} \cup A$ is 1 or greater for at least one set $A \subseteq \{v(1), \dots, v(i-1)\}$. (In our original version of the breadth-first ordering, $A$ would have been collapsed by now and so $w$ really is a child.)
6. Rename the children as before, remove them from $\mathcal{L}$ and add them to $\mathcal{R}$. Increment $i$ and repeat from step 4.

Observe that, before we find the children of $v(i)$, we do not look at the sets $\{v(i), w\} \cap A$, where $w \in \mathcal{L}$ and $A \subseteq \{v(1), \dots, v(i-1)\}$. Thus, in order to find the children, we need only consider random variables which are independent of what we have seen so far. So, if we imagine deleting $v(1), \dots, v(i-1)$, the 2-edge parameter in the remaining hypergraph is, indeed, $\lambda_2(N, i-1)$. Thus, any particular 2-edge is present with probability $\rho(N, i-1) = 1 - \exp(-\lambda_2(N, i-1))$.

Finally, we need to find the number of vertices eligible to be children of $v(i)$:

$$\begin{aligned}
&\#\{\text{vertices which cannot be a child of } v(i)\} \\
&\quad = |\mathcal{R}| \\
&\quad = C_N(1) + \cdots + C_N(i-1) + P_N(i-1) \\
&\quad = Z_N(i-1) + (i-1) + P_N(i-1),
\end{aligned}$$

and so $C_N(i) \sim \mathrm{Bin}(N - (i-1) - Z_N(i-1) - P_N(i-1), \rho(N, i-1))$. $\square$

In the statement of Theorem 3.1, we have used the floor function to interpolate between integer-valued time-points. Here, we will prove that the process $(Z_N(i))_{0 \le i \le N}$ converges with a different interpolation but this will be equivalent to the theorem as stated. Let $(E_{i,j} : 1 \le i \le N,\ 1 \le j \le N - (i-1) - Z_N(i-1) - P_N(i-1))$ be a family of independent $\mathrm{Exp}(\lambda_2(N, i-1))$

random variables and set

$$Z_N(i-1+u) = Z_N(i-1) - u + \sum_{j=1}^{N-(i-1)-Z_N(i-1)-P_N(i-1)} \mathbb{1}_{\{E_{i,j}\le u\}}, \qquad 0 \le u \le 1.$$

Take the filtration $\mathcal{F}_t^N = \sigma(Z_N(u) : u \le t)$. This filtration is spatial in that it tells us what we can "see" on the hypergraph at time $t$. Imagine that the vertex $v(i)$ is deleted at time $i-1$ [recall that $Z_N(i)$ is the number of patches on the hypergraph *after* the deletion of $v(i)$, adjusted for patches we have introduced]. Imagine that the patches on the children of $i$ appear one by one in the interval $(i-1,1]$. There are $N-(i-1)-Z_N(i-1)-P_N(i-1)$ possible children of $v(i)$ (i.e., vertices we have not yet reached) and each of them is actually a child with probability $\rho(N,i-1) = 1-\exp(-\lambda_2(N,i-1))$. Imagine each of these potential children having an $\operatorname{Exp}(\lambda_2(N,i-1))$ random variable associated with it; then the ones which actually are children are those whose exponential random variable is less than 1.

Define new processes by the standard decompositions

$$Z_N = M_N + A_N, \tag{6.1}$$

$$M_N^2 = R_N + Q_N, \tag{6.2}$$

where $M_N(0) = A_N(0) = R_N(0) = Q_N(0) = 0$, $M_N$ and $R_N$ are martingales, $A_N$ is a continuous adapted bounded variation process and $Q_N$ is a continuous adapted increasing process.

We will show that, for fixed $t_0$,

$$\frac{1}{N^{\alpha(k)/2}} \sup_{t\le N^{\alpha(k)}t_0} \left| A_N(t) - \frac{\mu_k}{(k-1)} \frac{t^{k-1}}{N^{k-2}} \right| \xrightarrow{p} 0, \tag{6.3}$$

$$\frac{1}{N^{\alpha(k)}} Q_N(N^{\alpha(k)}t_0) \xrightarrow{p} t_0, \tag{6.4}$$

$$\frac{1}{N^{\alpha(k)}} \mathbb{E}\left[ \sup_{t\le N^{\alpha(k)}t_0} |M_N(t) - M_N(t-)|^2 \right] \to 0. \tag{6.5}$$

Equivalently, on rescaling, we will show that

$$\sup_{t\le t_0} \left| \bar{A}_N(t) - \frac{\mu_k}{k-1} t^{k-1} \right| \xrightarrow{p} 0,$$

$$\bar{Q}_N(t_0) \xrightarrow{p} t_0,$$

$$\mathbb{E}\left[ \sup_{t\le t_0} |\bar{M}_N(t) - \bar{M}_N(t-)|^2 \right] \to 0,$$

where

$$\bar{A}_N(t) = N^{-\alpha(k)/2} A_N(N^{\alpha(k)} t),$$
$$\bar{M}_N(t) = N^{-\alpha(k)/2} M_N(N^{\alpha(k)} t),$$
$$\bar{Q}_N(t) = N^{-\alpha(k)} Q_N(N^{\alpha(k)} t).$$

By the martingale central limit theorem (Theorem 7.1.4(b) of [4]), conditions (6.4) and (6.5) are sufficient to prove that

$$\bar{M}_N \xrightarrow{d} B,$$

where $B$ is a standard Brownian motion. In conjunction with (6.3), this implies that

$$(\bar{Z}_N(t))_{t \geq 0} \xrightarrow{d} \left( B(t) + \frac{\mu_k}{k-1} t^{k-1} \right)_{t \geq 0}.$$

Now, since $A_N$ is continuous, the jumps of $M_N$ are those of $Z_N$, which are of size 1 by construction and so (6.5) is obvious. So we just need to find the explicit forms of the martingale decompositions.

Define $P_N(t) = P_N(\lfloor t \rfloor) = 1 - \inf_{s \leq \lfloor t \rfloor} Z_N(s)$, a continuous-time version of number of patches added to keep the process going. Let

$$a_N(t)\, dt = \mathbb{P}(\text{some new child appears during } [t, t+dt] | \mathcal{F}_t^N).$$

Then, as $Z_N$ has drift of rate $-1$ and jumps of $+1$,

$$A_N(t) = \int_0^t (a_N(s) - 1)\, ds, \tag{6.6}$$

$$Q_N(t) = \int_0^t a_N(s)\, ds. \tag{6.7}$$

Heuristically, this is because

$$\begin{aligned} &A_N(t+dt) - A_N(t) \\ &\qquad = \mathbb{E}[Z_N(t+dt) - Z_N(t) | \mathcal{F}_t^N] \\ &\qquad = -dt + \mathbb{P}(\text{some new child appears during } [t, t+dt] | \mathcal{F}_t^N), \\ &Q_N(t+dt) - Q_N(t) \\ &\qquad = \mathbb{E}[(Z_N(t+dt) - Z_N(t))^2 | \mathcal{F}_t^N] \\ &\qquad = \mathbb{P}(\text{some new child appears during } [t, t+dt] | \mathcal{F}_t^N). \end{aligned}$$

Thus,

$$Q_N(t) = A_N(t) + t. \tag{6.8}$$

Define $\lambda_2(N,t) = \lambda_2(N,\lfloor t \rfloor)$. So, for $0 \le u \le 1$,

$$\begin{aligned}
&\#\{\text{vertices at time } i-1+u \text{ which cannot be a child of } v(i)\}\\
&\quad = Z_N(i-1) + (i-1) + P_N(i-1) + \sum_{j=1}^{N-(i-1)-Z_N(i-1)-P_N(i-1)} \mathbb{1}_{\{E_{i,j}\le u\}}\\
&\quad = Z_N(i-1+u) + (i-1+u) + P_N(i-1+u),
\end{aligned}$$

and so, at time $t$, the time until the next child appears is the minimum of $N-t-Z_N(t)-P_N(t)$ independent $\mathrm{Exp}(\lambda_2(N,t))$ random variables, and so is exponential with parameter $(N-t-Z_N(t)-P_N(t))\lambda_2(N,t)$. Hence,

$$a_N(t) = (N-t-Z_N(t)-P_N(t))\lambda_2(N,t).$$

Then we have the following lemma:

Lemma 6.2. *Fix $0 \le r < 1$. Then for $N$ sufficiently large and all $0 \le t \le rN$,*

$$\begin{aligned}
&\left| N\lambda_2(N,t) - \left( \sum_{i=0}^{k-3} \left(\frac{t}{N}\right)^i + k(k-1)\beta_k \left(\frac{t}{N}\right)^{k-2} \right) \right| \\
&\qquad = \mathcal{O}((\log N)^2/N + (t/N)^{k-1}),
\end{aligned} \tag{6.9}$$

$$\begin{aligned}
&\left| a_N(t) - 1 - \mu_k \left(\frac{t}{N}\right)^{k-2} \right| \\
&\qquad \le \beta(r) \frac{|Z_N(t)| + P_N(t)}{N} + \mathcal{O}\left( \frac{(\log N)^2}{N} + \left(\frac{t}{N}\right)^{k-1} \right).
\end{aligned} \tag{6.10}$$

Proof. Lemma 6.1 of [3] says that there exists a constant $C < \infty$ such that

$$|N\lambda_2(N,t) - \beta''(t/N)| \le C(\log N)^2/N \tag{6.11}$$

for all $N \in \mathbb{N}$ and all $t \in \mathbb{Z} \cap [0, rN]$. Also, for $0 \le t \le rN$,

$$\begin{aligned}
\left| \beta''(t/N) - \left( \sum_{i=0}^{k-3} (t/N)^i + k(k-1)\beta_k (t/N)^{k-2} \right) \right| &\le (t/N)^{k-1} \sup_{0\le s<r} \beta^{(k+1)}(s)\\
&= (t/N)^{k-1}\beta^{(k+1)}(r)
\end{aligned}$$

by Taylor's theorem and (6.9) follows. Furthermore, using (6.9),

$$a_N(t) = \left(1 - \frac{t}{N} - \frac{Z_N(t)+P_N(t)}{N}\right) N\lambda_2(N,t)$$

$$= \left(1 - \frac{t}{N} - \frac{Z_N(t) + P_N(t)}{N}\right)\left(\sum_{i=0}^{k-3}\left(\frac{t}{N}\right)^i + k(k-1)\beta_k\left(\frac{t}{N}\right)^{k-2}\right)$$
$$+ \mathcal{O}\left(\frac{(\log N)^2}{N} + \left(\frac{t}{N}\right)^{k-1}\right)$$

and so

$$\left|a_N(t) - 1 - \mu_k\left(\frac{t}{N}\right)^{k-2}\right| \le \beta(r)\frac{|Z_N(t)| + P_N(t)}{N} + \mathcal{O}\left(\frac{(\log N)^2}{N} + \left(\frac{t}{N}\right)^{k-1}\right)$$

for sufficiently large $N$ and (6.10) holds. □

Thus, on integrating $a_N(t) - 1 - \mu_k(t/N)^{k-2}$ and using the fact that $P_N(t) = 1 - \inf_{s \le \lfloor t \rfloor} Z_N(s)$, we obtain that for some constant $C < \infty$,

$$\left|A_N(t) - \frac{\mu_k}{(k-1)}\frac{t^{k-1}}{N^{k-2}}\right| \le \frac{2Ct \max_{s \le t}|Z_N(s)|}{N} + \mathcal{O}\left(\frac{t(\log N)^2}{N} + \frac{t^k}{N^{k-1}}\right).$$

We wish to prove that

$$\frac{1}{N^{\alpha(k)/2}} \sup_{t \le N^{\alpha(k)} t_0}\left|A_N(t) - \frac{\mu_k}{(k-1)}\frac{t^{k-1}}{N^{k-2}}\right| \xrightarrow{p} 0$$

and so it will be sufficient to prove that

$$\frac{1}{N^{1-\alpha(k)/2}} \sup_{t \le N^{\alpha(k)} t_0} |Z_N(t)| \xrightarrow{p} 0$$

or the stronger statement that $\frac{1}{N^{\alpha(k)/2}} \sup_{t \le N^{\alpha(k)} t_0} |Z_N(t)|$ is stochastically bounded as $N \to \infty$.

Fix a large constant $K$ and let $T_N = \inf\{t \ge 0 : |Z_N(t)| > KN^{\alpha(k)/2}\} \wedge (N^{\alpha(k)} t_0)$. Then

$$\begin{aligned}
\mathbb{E}|Z_N(T_N)| &\le \mathbb{E}|M_N(T_N)| + \mathbb{E}|A_N(T_N)| \\
&\le \sqrt{\mathbb{E}[M_N^2(T_N)]} + \mathbb{E}|A_N(T_N)| \\
&= \sqrt{\mathbb{E}[Q_N(T_N)]} + \mathbb{E}|A_N(T_N)| \\
&= \sqrt{\mathbb{E}[A_N(T_N) + T_N]} + \mathbb{E}|A_N(T_N)| \\
&\le \sqrt{\mathbb{E}|A_N(T_N)| + N^{\alpha(k)} t_0} + \mathbb{E}|A_N(T_N)|,
\end{aligned}$$

where the equality on the third line is by the optional stopping theorem applied to (6.2) and the fourth line is from (6.8). Now, there exists a constant

$C < \infty$ such that

$$\begin{aligned}
\mathbb{E}|A_N(T_N)| &\le \mathbb{E}\left[\int_0^{T_N} |a_N(t) - 1|\, dt\right] \\
&\le \mathbb{E}\left[\int_0^{N^{\alpha(k)}t_0} \mu_k \left(\frac{t}{N}\right)^{k-2} dt\right] + \mathbb{E}\left[\int_0^{T_N} C\frac{P_N(t) + |Z_N(t)|}{N}\, dt\right] \\
&\quad + \mathcal{O}(N^{(k-3)/(2k-3)}) \\
&\le \mathbb{E}\left[\frac{2CT_N}{N} \sup_{t<T_N} |Z_N(t)|\right] + \mathcal{O}(N^{\alpha(k)/2}) \\
&\le \frac{2CN^{\alpha(k)}t_0 K N^{\alpha(k)/2}}{N} + \mathcal{O}(N^{\alpha(k)/2}) \\
&= 2Ct_0 K N^{(k-3)/(2k-3)} + \mathcal{O}(N^{\alpha(k)/2}).
\end{aligned}$$

Hence, $\mathbb{E}|Z_N(T_N)| \le C'KN^{(k-3)/(2k-3)} + C'N^{\alpha(k)/2}$ for some constant $C'$ and so, by Markov's inequality,

$$\begin{aligned}
\mathbb{P}\left(N^{-\alpha(k)/2} \sup_{t \le N^{\alpha(k)}t_0} |Z_N(t)| > K\right) &= \mathbb{P}(|Z_N(T_N)| > KN^{\alpha(k)/2}) \\
&\le \frac{\mathbb{E}|Z_N(T_N)|}{KN^{\alpha(k)/2}} \\
&\le \frac{C'}{K} + C'N^{-1/(2k-3)},
\end{aligned}$$

which gives the required stochastic boundedness as $N \to \infty$.

Finally, we need to show that

$$\frac{1}{N^{\alpha(k)}} Q_N(N^{\alpha(k)}t_0) \xrightarrow{p} t_0.$$

But since $Q_N(t) = A_N(t) + t$, this follows immediately from

$$\frac{1}{N^{\alpha(k)/2}} \sup_{t \le N^{\alpha(k)}t_0} \left| A_N(t) - \frac{\mu_k}{(k-1)} \frac{t^{k-1}}{N^{k-2}} \right| \xrightarrow{p} 0. \qquad \square$$

Proof of Theorem 5.2. Return now to the discrete-time setting of Theorem 5.2. The proof is quite involved and so we will begin by outlining our method. The random walk $Z_N$, started from above $RN^{\alpha(k)/2}$ at time $R^2N^{\alpha(k)}$, is close to dominating the breadth-first walk associated with the tree produced by a Galton–Watson branching process with Poisson offspring distribution of mean $1 + \frac{1}{2}\mu_k N^{-\alpha(k)/2} R^{2(k-2)}$. This branching process is ( just) supercritical and so the associated breadth-first walk has both positive drift and the virtue of identically distributed step sizes. We will essentially work with that process instead, using it as a lower bound for our original breadth-first walk. If, with high probability, the lower bound goes for $N\delta$ steps without hitting 0, then so does our original breadth-first walk. In order to show that the smaller random walk is unlikely to come back to 0, we use an exponential supermartingale argument and the optional stopping theorem.

Fix $R > 0$ and work on the set

$$\Omega_{N,R} = \{Z_N(\lceil N^{\alpha(k)} R^2 \rceil) > RN^{\alpha(k)/2}, P_N(\lceil N^{\alpha(k)} R^2 \rceil) < RN^{\alpha(k)/2}\}.$$

As we have already seen at (5.1) and (5.2), this set has asymptotically high probability. It will be useful to keep track of the times when $(Z_N(i))_{\lceil N^{\alpha(k)} R^2 \rceil \le i \le N}$ goes above and below the line $f(i) = \frac{1}{2}\mu_k i^{k-2} N^{-k+3}$. To this end, introduce two sequences of stopping times, $(\tau_n)_{n\ge1}$ and $(\tau'_n)_{n\ge1}$, such that

$$\tau_1 = \inf\{i \ge \lceil N^{\alpha(k)} R^2 \rceil : Z_N(i) \ge f(i)\}$$

and, for $n \ge 1$,

$$\tau'_n = \inf\{i > \tau_n : Z_N(i) < f(i)\},$$
$$\tau_{n+1} = \inf\{i > \tau'_n : Z_N(i) \ge f(i)\}.$$

Let $\Omega_N(i) = \{\lceil N^{\alpha(k)} R^2 \rceil \le i < \tau_1\} \cup (\bigcup_{n\ge1}\{\tau'_n \le i < \tau_{n+1}\})$, the event that $Z_N(i)$ is below the line. Define $\mathcal{F}_i^N = \sigma(Z_N(j) : j \le i)$.

LEMMA 6.3. *For sufficiently large $N$, there exists a $\theta < 0$ such that $L_N(i) = \exp(\theta Z_N(i))$ is a supermartingale on the disjoint sets of times*

$$\{\lceil N^{\alpha(k)} R^2 \rceil, \dots, \tau_1\}, \{\tau'_1, \dots, \tau_2\}, \{\tau'_2, \dots, \tau_3\}, \dots.$$

*That is, $(L_N(i))_{0\le i\le N}$ is an integrable adapted process and*

$$\mathbb{E}[L_N(i) | \mathcal{F}_{i-1}^N] \le L_N(i-1)$$

*whenever $i - 1 \in \{\lceil N^{\alpha(k)} R^2 \rceil, \dots, \tau_1 - 1\} \cup (\bigcup_{n\ge1}\{\tau'_n, \dots, \tau_{n+1} - 1\})$.*

PROOF. Consider the conditional moment generating function of an increment of $Z_N$: by Lemma 6.1,

$$\begin{aligned}
\phi_i(N,\theta) &= \mathbb{E}[\exp(\theta(Z_N(i) - Z_N(i-1))) | \mathcal{F}_{i-1}^N] \\
&= \exp\{(N - (i-1) - Z_N(i-1) - P_N(i-1)) \\
&\qquad \times \log(1 + \rho(N, i-1)(e^\theta - 1)) - \theta\} \\
&\le \exp\{(N - (i-1) - Z_N(i-1) - P_N(i-1)) \\
&\qquad \times \rho(N, i-1)(e^\theta - 1) - \theta\}.
\end{aligned}$$

Equation (6.11) implies that there exists a constant $C < \infty$ such that

$$N\rho(N, i-1) \geq \beta''\left(\frac{i}{N}\right) - \frac{C(\log N)^2}{N}.$$

On the set $\Omega_{N,R} \cap \Omega_N(i-1)$, we have $P_N(i-1) < RN^{\alpha(k)/2}$ and so, for sufficiently large $N$,

$$\begin{aligned}&(N-(i-1)-Z_N(i-1)-P_N(i-1))\rho(N,i-1)\\&\quad\geq\left(1-\frac{i}{N}-\frac{1}{2}\mu_k\left(\frac{i}{N}\right)^{k-2}\right)\beta''\left(\frac{i}{N}\right)-KN^{-(k-1)/(2k-3)}\end{aligned}$$

for some constant $K < \infty$ depending on $R$. The first nonzero derivative of $(1-s-\frac{1}{2}\mu_k s^{k-2})\beta''(s)$ at 0 is the $(k-2)$nd which is $\frac{1}{2}(k-2)!\mu_k$. Thus, there exists a $\delta > 0$ such that $(1-s-\frac{1}{2}\mu_k s^{k-2})\beta''(s)$ is increasing on $[0,\delta]$. Therefore, for $i \leq \delta N$ on $\Omega_{N,R} \cap \Omega_N(i-1)$, we have

$$\begin{aligned}&\left(1-\frac{i}{N}-\frac{1}{2}\mu_k\left(\frac{i}{N}\right)^{k-2}\right)\beta''\left(\frac{i}{N}\right)\\&\quad\geq\left(1-N^{\alpha(k)-1}R^2-\frac{1}{2}\mu_k N^{(\alpha(k)-1)(k-2)}R^{2(k-2)}\right)\beta''(N^{\alpha(k)-1}R^2)\end{aligned}$$

by putting in $i = N^{\alpha(k)}R^2$. Now, expanding $\beta''(s)$ as $1+s+\cdots+s^{k-3}+k(k-1)\beta_k s^{k-2}$ and incurring an error of size $\mathcal{O}(s^{k-1})$, we see that the right-hand side is bounded below by

$$1+\tfrac{1}{2}\mu_k N^{-\alpha(k)/2}R^{2(k-2)} - K'N^{-(k-1)/(2k-3)}$$

for some constant $K' < \infty$. Thus,

$$\begin{aligned}\phi_i(N,\theta) \leq \exp\{&(1+\tfrac{1}{2}\mu_k N^{-\alpha(k)/2}R^{2(k-2)})(e^\theta-1)-\theta\\&+K'N^{-(k-1)/(2k-3)}(1-e^\theta)\}.\end{aligned}$$

If we had the breadth-first walk on the tree produced by a branching process with Poisson offspring distribution of mean $1+\frac{1}{2}\mu_k N^{-\alpha(k)/2}R^{2(k-2)}$, we would have

$$\exp\{(1+\tfrac{1}{2}\mu_k N^{-\alpha(k)/2}R^{2(k-2)})(e^\theta-1)-\theta\}$$

as the conditional moment generating function of an increment. Thus, we will effectively use this simpler process as a "lower bound" for our original process.

Now, let $\bar{\theta}$ be the value of $\theta$ which minimizes

$$\exp\{(1+\tfrac{1}{2}\mu_k N^{-\alpha(k)/2}R^{2(k-2)})(e^\theta-1)-\theta\},$$

so that it is easily seen that

$$\bar{\theta} = -\log(1+\tfrac{1}{2}\mu_k N^{-\alpha(k)/2}R^{2(k-2)})$$

(and so, trivially, $\bar{\theta} < 0$). For $i \le N\delta$ on $\Omega_{N,R} \cap \Omega_N(i-1)$,

$$\begin{aligned}\phi_i(N,\bar{\theta}) &\le \exp\biggl\{\biggl(1+\frac{1}{2}\mu_k N^{-\alpha(k)/2} R^{2(k-2)}\biggr)(e^{\bar{\theta}}-1)-\bar{\theta} \\ &\qquad\qquad + K' N^{-(k-1)/(2k-3)}(1-e^{\bar{\theta}})\biggr\} \\ &= \exp\biggl\{\log\biggl(1+\frac{1}{2}\mu_k N^{-\alpha(k)/2}R^{2(k-2)}\biggr) - \frac{1}{2}\mu_k N^{-\alpha(k)/2} R^{2(k-2)} \\ &\qquad\qquad + \frac{(1/2)\mu_k K' N^{-1} R^{2(k-2)}}{1+(1/2)\mu_k N^{-\alpha(k)/2} R^{2(k-2)}}\biggr\} \\ &\le \exp\{-C_2 N^{-\alpha(k)} + C_1 N^{-1}\} \\ &\le 1,\end{aligned}$$

for some constants $C_1, C_2$ and sufficiently large $N$. Hence, for $\theta = \bar{\theta}$ and sufficiently large $N$, $\mathbb{E}[L_N(i)|\mathcal{F}_{i-1}^N] \le L_N(i-1)$ on the set $\Omega_{N,R} \cap \Omega_N(i-1)$. □

LEMMA 6.4. *We have*

$$\lim_{N\to\infty} \mathbb{P}(S_N^R \le N\delta | \Omega_{N,R}) \le \exp(-\tfrac{1}{2}\mu_k R^{2k-1}).$$

PROOF. $Z_N$ cannot hit 0 when it is above the line $f(i) = \frac{1}{2}\mu_k i^{k-2} N^{-k+3}$. Thus, if $Z_N$ does hit 0 before time $N\delta$, it must occur before time $\tau_1$ or between times $\tau_i'$ and $\tau_{i+1}$ for some $i \ge 1$ and so

$$\begin{aligned}&\mathbb{P}(S_N^R \le N\delta|\Omega_{N,R}) \\ &\quad\le \mathbb{P}(S_N^R \le \tau_1 \wedge (N\delta)|\Omega_{N,R}) \\ &\qquad + \sum_{i=1}^{\infty} \mathbb{P}(S_N^R \le \tau_{i+1} \wedge (N\delta)|\Omega_{N,R}, \tau_i < S_N^R \wedge (N\delta), \tau_i' < N\delta) \\ &\quad\le \mathbb{P}(S_N^R \le \tau_1 \wedge (N\delta)|\Omega_{N,R}) \\ &\qquad + \sum_{i=\lceil N^{\alpha(k)}R^2\rceil}^{\infty} \mathbb{P}(S_N^R \le T_i \wedge (N\delta)|\Omega_{N,R}, Z_N(i-1) \ge f(i-1), \\ &\qquad\qquad\qquad\qquad Z_N(i) < f(i)),\end{aligned}$$

where $T_i = \inf\{j \ge i : Z_N(j) > f(j)\}$ and each term in the above summation expresses the probability of a downcrossing from the line $f(i)$ to 0. It will turn out that only the term $\mathbb{P}(S_N^R \le \tau_1 \wedge (N\delta)|\Omega_{N,R})$ makes a significant contribution (intuitively because the process is much closer to 0 at time $\lceil N^{\alpha(k)} R^2 \rceil$ than it is at $\tau_i'$ for any $i \ge 1$).

As $S_N^R \wedge \tau_1 \wedge (N\delta)$ is a bounded stopping time, by the optional stopping theorem we obtain that

$$\begin{aligned} L_N(\lceil N^{\alpha(k)}R^2\rceil)\mathbb{1}_{\Omega_{N,R}} &\geq \mathbb{E}[L_N(S_N^R \wedge \tau_1 \wedge (N\delta))\mathbb{1}_{\Omega_{N,R}}|\mathcal{F}^N_{\lceil N^{\alpha(k)}R^2\rceil}] \\ &\geq \mathbb{P}(S_N^R \leq \tau_1 \wedge (N\delta)|\mathcal{F}^N_{\lceil N^{\alpha(k)}R^2\rceil})\mathbb{1}_{\Omega_{N,R}}. \end{aligned}$$

Hence,

$$\begin{aligned} \mathbb{P}(S_N^R \leq \tau_1 \wedge (N\delta)|\Omega_{N,R}) &\leq \mathbb{E}[L_N(\lceil N^{\alpha(k)}R^2\rceil)|\Omega_{N,R}] \\ &= \mathbb{E}[\exp(\bar{\theta}Z_N(\lceil N^{\alpha(k)}R^2\rceil))|\Omega_{N,R}] \\ &\leq \exp(\bar{\theta}RN^{\alpha(k)/2}) \\ &= \exp(-\log(1+\tfrac{1}{2}\mu_k N^{-\alpha(k)/2}R^{2(k-2)})RN^{\alpha(k)/2}) \\ &= \exp(-\tfrac{1}{2}\mu_k R^{2k-1} + \mathcal{O}(N^{-\alpha(k)/2})) \\ &\to \exp(-\tfrac{1}{2}\mu_k R^{2k-1}) \end{aligned}$$

as $N \to \infty$. By a similar argument,

$$\begin{aligned} &\mathbb{P}(S_N^R \leq T_i \wedge (N\delta)|\Omega_{N,R}, Z_N(i-1) \geq f(i-1), Z_N(i) < f(i)) \\ &\quad \leq \mathbb{E}[L_N(i)|\Omega_{N,R}, Z_N(i-1) \geq f(i-1), Z_N(i) < f(i)] \\ &\quad \leq \mathbb{E}[L_N(i)|\Omega_{N,R}, Z_N(i) = \lceil f(i-1)\rceil - 1] \\ &\quad \leq C\exp(\tfrac{1}{2}\mu_k\bar{\theta}i^{k-2}N^{-k+3}), \end{aligned}$$

where $C$ is a constant and the second inequality holds because $Z_N$ can step down by at most 1 and so the smallest that $Z_N(i)$ can be and still have had $Z_N(i-1)$ above the line is $\lceil f(i-1)\rceil - 1$. For $i \geq \lceil N^{\alpha(k)}R^2\rceil$,

$$\begin{aligned} &\exp(\tfrac{1}{2}\mu_k\bar{\theta}i^{k-2}N^{-k+3}) \\ &\quad \leq \exp(\tfrac{1}{2}\mu_k\bar{\theta}N^{-(k-3)/(2k-3)}R^{2(k-3)}i) \\ &\quad = \exp(-\tfrac{1}{2}\mu_k N^{-(k-3)/(2k-3)}R^{2(k-3)}\log(1+\tfrac{1}{2}\mu_k N^{-\alpha(k)/2}R^{2(k-3)})i) \\ &\quad = [(1+\tfrac{1}{2}\mu_k N^{-\alpha(k)/2}R^{2(k-3)})^{-(1/2)\mu_k N^{-(k-3)/(2k-3)}R^{2(k-3)}}]^i. \end{aligned}$$

Let

$$g(N,R) = (1+\tfrac{1}{2}\mu_k N^{-\alpha(k)/2}R^{2(k-3)})^{-(1/2)\mu_k N^{-(k-3)/(2k-3)}R^{2(k-3)}}.$$

Then

$$\begin{aligned} &\sum_{i=\lceil N^{\alpha(k)}R^2\rceil}^{\infty} \mathbb{P}(S_N^R \leq T_i \wedge (N\delta)|\Omega_{N,R}, Z_N(i-1) \geq f(i-1), Z_N(i) < f(i)) \\ &\quad \leq C\sum_{i=\lceil N^{\alpha(k)}R^2\rceil}^{\infty} g(N,R)^i = \frac{Cg(N,R)^{\lceil N^{\alpha(k)}R^2\rceil}}{1-g(N,R)}. \end{aligned}$$

This behaves essentially like $N^{(2k-5)/(2k-3)}\exp(-N^{1/(2k-3)})$ and so converges to 0 as $N \to \infty$. Hence,

$$\lim_{N\to\infty} \mathbb{P}(S_N^R \le N\delta | \Omega_{N,R}) \le \exp(-\tfrac{1}{2}\mu_k R^{2k-1}). \qquad \square$$

Now note that

$$\mathbb{P}(S_N^R \le N\delta) \le \mathbb{P}(S_N^R \le N\delta | \Omega_{N,R}) + \mathbb{P}(\Omega_{N,R}^c)$$

and so

$$\lim_{N\to\infty} \mathbb{P}(S_N^R \le N\delta) \le \exp(-\tfrac{1}{2}\mu_k R^{2k-1}) + \mathbb{P}(W^k(R^2) \le R) + \mathbb{P}\Big(-\inf_{s\ge 0} W^k(s) \ge R\Big),$$

which converges to 0 as $R \to \infty$. Theorem 5.2 follows. $\square$

Proof of Theorem 5.3. We will make use of a weaker version of the fluid limit methods expounded in [3]. Suppose, for the moment, that $(X_t^N)_{t\ge 0}$ is a time-homogeneous pure jump Markov process taking values in $I^N = \frac{1}{N}\mathbb{Z}^d \subseteq \mathbb{R}^d$. Let $K^N(x,dy)$ be the Lévy kernel of $(X_t^N)_{t\ge 0}$. Define the Laplace transform of this Lévy kernel by

$$m^N(x,\theta) = \int_{\mathbb{R}^d} e^{\langle \theta, y\rangle} K^N(x,dy),$$

where $\langle \cdot,\cdot\rangle$ is the usual inner product on $\mathbb{R}^d$. Let $S$ be an open subset of $\mathbb{R}^d$ and define $S^N = I^N \cap S$. We are interested in the behavior of $(X_t^N)_{t\ge 0}$ up until the first time it leaves $S$ (e.g., this may mean that one of its coordinates hits 0). With this in mind, define the relevant stopping time

$$T^N = \inf\{t \ge 0 : X_t^N \notin S\} \wedge t_0,$$

where $t_0 > 0$ is a constant.

Theorem 6.5 ([3]). *Assume the following conditions:*

1. *There exists a limit kernel $K(x,dy)$, defined for $x \in S$, and a constant $\eta_0 > 0$ such that $m(x,\theta) < \infty$ for all $x \in S$ and $\|\theta\| \le \eta_0$, where $m$ is the Laplace transform of $K$.*
2. *We have*

$$\sup_{x\in S^N} \sup_{\|\theta\|\le\eta_0} |N^{-1} m^N(x, N\theta) - m(x,\theta)| \to 0$$

   *as $N \to \infty$.*
3. *Let $b(x) = m'(x,0)$, where $m'(x,\theta)$ is the vector of partial derivatives in components of $\theta$. Then assume that $b$ is Lipschitz on $S$.*

4. *We have* $\|X_0^N - x_0\| \xrightarrow{p} 0$ *as* $N \to \infty$ *for some constant* $x_0 \in \bar{S}$.

*Denote by* $(x(t))_{t\geq 0}$ *the unique solution to the differential equation* $\dot{x}(t) = \tilde{b}(x(t))$ *with initial condition* $x(0) = x_0$, *where* $\tilde{b}$ *is a Lipschitz extension of* $b$ *to* $\mathbb{R}^d$. *Then,*

$$\sup_{0\leq t\leq T^N} \|X_t^N - x(t)\| \xrightarrow{p} 0$$

*as* $N \to \infty$.

In simple cases, where $x(t)$ does not graze the boundary of $S$ before crossing it, it is straightforward to show that $T^N$ converges in probability to the first exit time of $x(t)$ from $S$.

In order to apply this theorem, we need to be working with a pure jump Markov process. Now, $(Z_N(i), P_N(i))_{0\leq i\leq N}$ is a discrete-time Markov chain. Take $(\nu_t^N)_{t\geq 0}$ to be a Poisson process of rate $N$ and let

$$X_t^N = (X_t^{1,N}, X_t^{2,N}, X_t^{3,N}) = \frac{1}{N}(\nu_t^N, Z_N(\nu_t^N), P_N(\nu_t^N)).$$

We will prove a fluid limit with this embedding into continuous time, which will imply the theorem as stated.

We need to check that conditions 1–4 of Theorem 6.5 hold. The process $X^N$ is naturally defined on $I^N = \{x \in \mathbb{R}^3 : Nx_1 \in \mathbb{Z}^+, Nx_2 \in \mathbb{Z}, Nx_3 \in \mathbb{Z}^+\}$. Let $S = \{x \in R^3 : |x_1| < r_1, 0 < x_2 < r_2, |x_3| < r_3\}$ for constants $r_1, r_2, r_3 < \infty$ and let $S^N = I^N \cap S$. Let $K^N(x, dy)$ be the Lévy kernel of $X^N$. Then, using the representation (3.1) of the evolution of $P_N$ and Lemma 6.1, $N^{-1}K^N(x, \cdot)$ is the law of

$$N^{-1}(1, B_N - 1, \mathbb{1}_{\{x_2+x_3=1/N,\ B_N=0\}}),$$

where $B_N \sim \mathrm{Bin}(N - Nx_1 - Nx_2 - Nx_3, \rho(N, Nx_1))$. Thus, $K^N(x, dy)$ has Laplace transform

$$\begin{aligned} m_1^N(x,\theta) &= Ne^{\theta/N}, \\ m_2^N(x,\theta) &= N\exp((N - Nx_1 - Nx_2 - Nx_3) \\ &\qquad\qquad \times \log(1 + \rho(N, Nx_1)(e^{\theta/N} - 1)) - \theta/N), \\ m_3^N(x,\theta) &= \begin{cases} Ne^{\theta/N}(1-\rho(N,Nx_1))^{N-Nx_1-1} \\ \quad + N(1-(1-\rho(N,Nx_1))^{N-Nx_1-1}), & \text{if } x_2 + x_3 = 1/N, \\ N, & \text{if } x_2 + x_3 \neq 1/N. \end{cases} \end{aligned}$$

Using (6.11), $B_N \xrightarrow{d} \mathrm{Poisson}((1 - x_1 - x_2 - x_3)\beta''(x_1))$ as $N \to \infty$ and so there exists a limit kernel $K(x, dy)$ with Laplace transform

$$m_1(x,\theta) = e^{\theta},$$

$$m_2(x,\theta) = \exp((1 - x_1 - x_2 - x_3)\beta''(x_1)(e^\theta - 1) - \theta),$$

$$m_3(x,\theta) = \begin{cases} (e^\theta - 1)\exp(-(1 - x_1)\beta''(x_1)) + 1, & \text{if } x_2 = -x_3, \\ 1, & \text{if } x_2 \neq -x_3. \end{cases}$$

Furthermore, there exists $\eta_0 > 0$ such that $m(x,\theta) < \infty$ for all $x \in S$ and $\|\theta\| \leq \eta_0$ and also such that

$$\sup_{x\in S^N} \sup_{\|\theta\|<\eta_0} |N^{-1}m^N(x, N\theta) - m(x,\theta)| \to 0,$$

where $\|\cdot\|$ is the Euclidean norm on $\mathbb{R}^3$ (note that $x_2 \neq -x_3$ in $S$). Thus, conditions 1 and 2 are satisfied.

Let $T_N = \inf\{t \geq 0 : X_t^N \notin S^N\}$ and recall that $z(t) = 1 - t - \exp(-\beta'(t))$. Fix a large $R < \infty$. We will prove the fluid limit result in three time intervals: $[0, R^2 N^{\alpha(k)-1}]$, $[R^2 N^{\alpha(k)-1}, \delta]$ and $[\delta, 1]$.

*Time interval* $[0, R^2 N^{\alpha(k)-1}]$. Until time $R^2 N^{\alpha(k)-1}$, $Z_N$ is in the Brownian regime of Theorem 3.1 and so

$$\begin{aligned} &\sup_{0\leq t\leq R^2 N^{\alpha(k)-1}} |N^{-1}Z_N(\lfloor Nt\rfloor) - z(t)| \\ &\qquad \leq \sup_{0\leq t\leq R^2 N^{\alpha(k)-1}} |N^{-1}Z_N(Nt)| + \sup_{0\leq t\leq R^2 N^{\alpha(k)-1}} |z(t)|, \end{aligned}$$

which converges to 0 in probability as $N \to \infty$, regardless of the value of $R$. Similarly,

$$\sup_{0\leq t\leq R^2 N^{\alpha(k)-1}} |N^{-1}P_N(Nt)| \xrightarrow{p} 0.$$

It is elementary that $\sup_{0\leq t\leq R^2 N^{\alpha(k)-1}} |N^{-1}\nu_t^N - t| \xrightarrow{p} 0$. Thus,

$$\sup_{0\leq t\leq R^2 N^{\alpha(k)-1}} \|X_t^N - x(t)\| \xrightarrow{p} 0, \tag{6.12}$$

where $x(t) = (t, z(t), 0)$.

*Time interval* $[R^2 N^{\alpha(k)-1}, \delta]$. Expression (6.12) provides us with condition 4 of the fluid limit theorem. Suppose now that we fix $N$ and $R$ and work on the set

$$\Omega_{N,R} = \{X^{2,N}_{\lceil R^2 N^{\alpha(k)-1}\rceil} > RN^{\alpha(k)/2-1}, X^{3,N}_{\lceil R^2 N^{\alpha(k)-1}\rceil} < RN^{\alpha(k)/2-1}\}.$$

Define $T_N^R = \inf\{t \geq R^2 N^{\alpha(k)-1} : X_t^N \notin S\}$ (which, for $r_1$, $r_2$ and $r_3$ large enough in the definition of $S$, is the time that $X^{2,N}$ first hits 0; this is a Poissonized equivalent of $S_N^R$ of Section 5). On $\Omega_{N,R}$, $X_t^{3,N}$ is constant on the interval $[R^2N^{\alpha(k)-1}, T_N^R]$ and so

$$\sup_{0\le t\le T_N^R} |N^{-1}X_t^{3,N}| \xrightarrow{p} 0$$

also. Now, setting $b(x) = m'(x,0)$, we have that for $x \in S$,

$$\begin{aligned}
&b_1(x) = 1,\\
&b_2(x) = (1-x_1-x_2-x_3)\beta''(x_1) - 1,\\
&b_3(x) = 0.
\end{aligned}$$

Condition 3 is clearly satisfied. By taking $r_1, r_2$ and $r_3$ large enough, the differential equation $\dot{x}_t = b(x(t))$ has unique solution $x(t) = (t, 1-t-\exp(-\beta'(t)), 0) = (t, z(t), 0)$ in $S$. Thus, Theorem 6.5 entails that for all $\varepsilon > 0$,

$$\lim_{N\to\infty} \mathbb{P}\biggl(\sup_{R^2N^{\alpha(k)-1}\le t\le T_N^R\wedge\delta} \|X_t^N - x(t)\| > \varepsilon \bigg| \Omega_{N,R}\biggr) = 0.$$

Now, for small enough $\delta > 0$ we have that $\delta < t^*$ and $\lim_{R\to\infty}\lim_{N\to\infty}\mathbb{P}(S_N^R > N\delta) = 1$ and so we will also have $\lim_{R\to\infty}\lim_{N\to\infty}\mathbb{P}(T_N^R > \delta) = 1$. Then

$$\begin{aligned}
&\lim_{N\to\infty} \mathbb{P}\biggl(\sup_{R^2N^{\alpha(k)-1}\le t\le \delta} \|X_t^N - x(t)\| > \varepsilon\biggr)\\
&\quad\le \lim_{N\to\infty} \mathbb{P}\biggl(\sup_{R^2N^{\alpha(k)-1}\le t\le T_N^R\wedge\delta} \|X_t^N - x(t)\| > \varepsilon, T_N^R > \delta \bigg| \Omega_{N,R}\biggr)\\
&\qquad + \lim_{N\to\infty} \mathbb{P}(T_N^R \le \delta | \Omega_{N,R}) + \lim_{N\to\infty} \mathbb{P}(\Omega_{N,R}^c),
\end{aligned}$$

where the first of the three terms on the right-hand side is 0 as we have just proved.

Thus, for all $R < \infty$ and all $\varepsilon > 0$, we have that

$$\begin{aligned}
\lim_{N\to\infty} \mathbb{P}\biggl(\sup_{0\le t\le\delta} \|X_t^N - x(t)\| > \varepsilon\biggr) &\le \lim_{N\to\infty} \mathbb{P}\biggl(\sup_{0\le t\le R^2N^{\alpha(k)-1}} \|X_t^N - x(t)\| > \varepsilon\biggr)\\
&\quad + \lim_{N\to\infty} \mathbb{P}\biggl(\sup_{R^2N^{\alpha(k)-1}\le t\le\delta} \|X_t^N - x(t)\| > \varepsilon\biggr)
\end{aligned}$$

and so we can take the limit as $R \to \infty$ on the right-hand side to obtain

$$\begin{aligned}
&\lim_{N\to\infty} \mathbb{P}\biggl(\sup_{0\le t\le\delta} \|X_t^N - x(t)\| > \varepsilon\biggr)\\
&\quad\le \lim_{R\to\infty}\lim_{N\to\infty} \mathbb{P}(T_N^R \le \delta | \Omega_{N,R}) + \lim_{R\to\infty}\lim_{N\to\infty} \mathbb{P}(\Omega_{N,R}^c)\\
&\quad= 0.
\end{aligned}$$

*Time interval* $[\delta, 1]$. Suppose we now start from time $\delta$. We have just shown that the initial value converges in probability. Redefine $T_N = \inf\{t \geq \delta : X^N \notin S^N\}$ and let $\tilde{X}^N$ be such that $\tilde{X}_t^{1,N} = X_t^{1,N}$, $\tilde{X}_t^{2,N} = X_{t\wedge T_N}^{2,N}$ and $\tilde{X}_t^{3,N} = X_\delta^{3,N}$. Minor modifications to our argument in the previous section lead us to the limiting function $\tilde{x}(t) = (t, z(t \wedge t^*), 0)$ and the conclusion that for all $\varepsilon > 0$,

$$\lim_{N\to\infty} \mathbb{P}\left(\sup_{\delta \leq t \leq 1} \|\tilde{X}_t^N - \tilde{x}(t)\| > \varepsilon\right) = 0,$$

where we have that $T_N \xrightarrow{p} t^*$ by our assumption that there are no zeros of the function $z(t)$ in $(0, t^*)$. Thus, finally,

$$\begin{aligned}
&\lim_{N\to\infty} \mathbb{P}\left(\sup_{0\leq t\leq 1} |\tilde{X}_t^{2,N} - \tilde{z}(t)| > \varepsilon\right)\\
&\quad \leq \lim_{N\to\infty} \mathbb{P}\left(\sup_{0\leq t\leq \delta} |X_t^{2,N} - z(t)| > \varepsilon\right)\\
&\quad\quad + \lim_{N\to\infty} \mathbb{P}\left(\sup_{\delta\leq t\leq 1} |\tilde{X}_t^{2,N} - z(t\wedge t^*)| > \varepsilon\right)\\
&\quad = 0. \qquad \square
\end{aligned}$$

PROOF OF THEOREM 5.4. We use the interpolation from the proof of Theorem 3.1 and the Doob decompositions contained therein.

LEMMA 6.6. *Choose $0 < \sigma < t^*$. Then*

$$\frac{1}{\sqrt{N}}(Z_N(Nt) - A_N(Nt))_{0\leq t\leq\sigma} \xrightarrow{d} (G_t)_{0\leq t\leq \sigma}$$

*as $N \to \infty$, where $(G_t)_{t\geq 0}$ has the distribution of the time-changed standard Brownian motion $(B_{t+z(t)})_{t\geq 0}$.*

PROOF. By the martingale central limit theorem (Theorem 7.1.4(b) of [4]), it is sufficient to show that

$$\frac{1}{N} Q_N(Nt) \xrightarrow{p} z(t) + t$$

and that

$$\frac{1}{N}\mathbb{E}\left[\sup_{0\leq t\leq N\sigma} (M_N(t) - M_N(t-))^2\right] \to 0,$$

the latter being obvious from the fact that $M_N$ cannot jump by more than 1. Now, by (6.8), $Q_N(t) = A_N(t) + t$ and so it will be sufficient to show that

$\frac{1}{N}A_N(Nt) \xrightarrow{p} z(t)$. By (6.6) and the definition of $z$ at (5.3),

$$\begin{aligned}
&\frac{1}{N}A_N(Nt) - z(t)\\
&\quad = \frac{1}{N}\int_0^{Nt}((N-s-Z_N(s)-P_N(s))\lambda_2(N,s)-1)\,ds - \int_0^t \dot{z}(s)\,ds\\
&\quad = \int_0^t\bigg(\Big(1-s-\frac{1}{N}Z_N(Ns)-\frac{1}{N}P_N(Ns)\Big)N\lambda_2(N,Ns)\\
&\qquad\qquad\qquad - \beta''(s)(1-s-z(s))\bigg)\,ds,
\end{aligned}$$

for $t \le \sigma$. There exists a constant $C < \infty$ such that

$$\begin{aligned}
\left|\frac{1}{N}A_N(Nt) - z(t)\right| &\le \left|\int_0^t (1-s)(N\lambda_2(N,Ns)-\beta''(s))\,ds\right|\\
&\quad + C\int_0^t\left|\frac{1}{N}Z_N(Ns)-z(s)\right|ds + C\int_0^t\left|\frac{P_N(Ns)}{N}\right|ds.
\end{aligned}$$

By (6.11), there exists a $K < \infty$ such that the right-hand side is bounded above by

$$K\frac{(\log N)^2}{N} + C\int_0^t\left|\frac{1}{N}Z_N(Ns)-z(s)\right|ds + C\int_0^t\left|\frac{P_N(Ns)}{N}\right|ds.$$

As $\mu_k > 0$, by Theorems 3.1 and 5.2, only $\mathcal{O}(N^{\alpha(k)/2})$ patches need to be added before the breadth-first walk begins its giant excursion and so, as $\alpha(k) < 1$, we must have

$$\sup_{0\le t\le\sigma}\left|\frac{P_N(Nt)}{\sqrt{N}}\right| \xrightarrow{p} 0$$

as $N \to \infty$. Thus, $|N^{-1}A_N(Nt) - z(t)|$ converges to 0 in probability, uniformly in $0 \le t \le \sigma$, by Theorem 5.3. □

Now define

$$\begin{aligned}
X_t^N &= \frac{1}{\sqrt{N}}(Z_N(Nt) - Nz(t)),\\
G_t^N &= \frac{1}{\sqrt{N}}(Z_N(Nt) - A_N(Nt)).
\end{aligned}$$

Then we know already from Lemma 6.6 that $(G_t^N)_{0\le t\le\sigma} \xrightarrow{d} (G_t)_{0\le t\le\sigma}$ and we wish to prove that $(X_t^N)_{0\le t\le\sigma} \xrightarrow{d} (X_t)_{0\le t\le\sigma}$. Now,

$$\begin{aligned}
X_t^N &= G_t^N + \sqrt{N}\int_0^t (a_N(Ns)-1)\,ds - \sqrt{N}\int_0^t \dot{z}(s)\,ds\\
&= G_t^N + E_t^N - \int_0^t X_s^N\beta''(s)\,ds,
\end{aligned}$$

where

$$E_t^N = \int_0^t \left( \sqrt{N}(1-s) - \frac{1}{\sqrt{N}}(Z_N(Ns) + P_N(Ns)) \right) (N\lambda_2(N, Ns) - \beta''(s))\, ds$$
$$- \int_0^t \frac{1}{\sqrt{N}} P_N(Ns)\beta''(s)\, ds.$$

Observe that, because $1 \leq Z_N + P_N \leq N$, we have

$$|E_t^N| \leq \int_0^t 2\sqrt{N}|N\lambda_2(N, Ns) - \beta''(s)|\, ds + \beta''(\sigma)t \sup_{0\leq s\leq t}\left|\frac{P_N(Ns)}{\sqrt{N}}\right|$$
$$\leq 2t\frac{(\log N)^2}{\sqrt{N}} + \beta''(\sigma)t \sup_{0\leq s\leq t}\left|\frac{P_N(Ns)}{\sqrt{N}}\right|.$$

At the end of the proof of Lemma 6.6, we showed that $\sup_{0\leq t\leq\sigma}|\frac{P_N(Nt)}{\sqrt{N}}| \xrightarrow{p} 0$ as $N \to \infty$ and so $\sup_{0\leq t\leq\sigma}|E_t^N| \xrightarrow{p} 0$.

Now, there exists a continuous function (of *paths*) $F : C_0[0, \sigma] \to C_0[0, \sigma]$ such that $X^N = F(G^N + E^N)$. Let $X = F(G)$, that is, let $(X_t)$ satisfy

$$dX_t = dG_t - X_t\beta''(t)\, dt.$$

Then $(G_t^N + E_t^N)_{0\leq t\leq\sigma} \xrightarrow{d} (G_t)_{0\leq t\leq\sigma}$ and so, by continuity of $F$ and the continuous mapping theorem (Corollary 3.1.9 of [4]), $(X_t^N)_{0\leq t\leq\sigma} \xrightarrow{d} (X_t)_{0\leq t\leq\sigma}$. By using an integrating factor to solve the above stochastic differential equation, we see that

$$X_t = \exp(-\beta'(t)) \int_0^t \exp(\beta'(s))\, dG_s. \qquad \square$$

**Acknowledgments.** I am very grateful to James Norris for his support and encouragement and for many enlightening conversations. My thanks also to James Martin for some helpful discussions and to an anonymous referee for several suggestions leading to the improvement of the Introduction. The work for this paper was mostly done in the Statistical Laboratory at the University of Cambridge.

Laboratoire de Probabilités
et Modèles Aléatoires
Université Pierre et Marie Curie (Paris 6)
175 rue du Chevaleret
75013 Paris
France
e-mail: C.Goldschmidt@statslab.cam.ac.uk